\documentclass{gtart}

\def\ifplaintex{\expandafter\ifx\csname documentclass\endcsname\relax}

\def\gtp{{\mathsurround=0pt\it $\cal G\mskip-2mu$eometry \&\ 
$\cal T\!\!$opology $\cal P\!$ublications}}  

\def\recd{{\small Received:\qua\receiveddate\ifx\reviseddate\relax
\else\qquad Revised:\qua\reviseddate\fi\par}} 

\def\lognumber#1{\def\thelognumber{#1}}
\def\volumenumber#1{\def\thevolumenumber{#1}}
\def\volumeyear#1{\def\thevolumeyear{#1}}
\def\papernumber#1{\def\thepapernumber{#1}}
\def\pagenumbers#1#2{\def\startpage{#1}\def\finishpage{#2}}
\def\published#1{\def\publishdate{#1}}

\def\received#1{\def\receiveddate{#1}}

\def\accepted#1{\def\accepteddate{#1}}

\long\def\asciiabstract#1{\long\def\theasciiabstract{#1}}

\let\\\par\let\thelognumber\relax\let\thevolumenumber\relax
\let\thepapernumber\relax\let\thevolumeyear\relax\let\startpage\relax
\let\finishpage\relax\let\publishdate\relax\let\receiveddate\relax
\let\reviseddate\relax\let\accepteddate\relax\let\theasciititle\relax
\let\theasciiauthors\relax
\let\theasciiabstract\relax

\let\theasciiemail\relax

\ifplaintex
\font\logobig=cmssbx10 scaled 3836
\font\logomed=cmssbx10 scaled 2557
\else
\font\logobig=cmssbx10 scaled 4200
\font\logomed=cmssbx10 scaled 2800
\fi

\long\def\makeagttitle{   
\count0=\startpage
\agt\hfill      
\hbox to 45truept{\vbox to 0pt{\vglue -13truept{\logomed A\kern -.37em{\logobig 
T}\kern -.38em G}\vss}\hss}
\break
{\small Volume \thevolumenumber\ (\thevolumeyear)
\startpage--\finishpage\nl
Published: \publishdate}

\vglue .25truein

{\parskip=0pt\leftskip 0pt plus
1fil\def\\{\par\smallskip}{\Large\bf\thetitle}\par\medskip} \vglue
0.05truein

{\parskip=0pt\leftskip 0pt plus 1fil\def\\{\par}{\sc\theauthors}
\par\medskip}%
 
\vglue 0.03truein 

{\small\leftskip 25truept\rightskip 25truept{\bf Abstract}\stdspace\theabstract

{\bf AMS Classification}\stdspace\theprimaryclass
\ifx\thesecondaryclass\relax\else; \thesecondaryclass\fi\par
{\bf Keywords}\stdspace \thekeywords\par}\vglue 7truept

}   

\ifplaintex
\font\phead=cmsl9 scaled 950
\font\pnum=cmbx10 scaled 913
\font\pfoot=cmsl9 scaled 950
\headline{\vbox to 0pt{\vskip -4.5mm\line{\small\phead\ifnum
\count0=\startpage ISSN 1472-2739 (on-line) 1472-2747 (printed)
\hfill {\pnum\folio}\else\ifodd\count0\def\\{ }%
\ifx\theshorttitle\relax\thetitle\else\theshorttitle\fi\hfill{\pnum\folio}
\else\def\\{ and }{\pnum\folio}\hfill\ifx\theshortauthors\relax\theauthors
\else\theshortauthors\fi\fi\fi}\vss}}
\footline{\vbox to 0pt{\vglue 0mm\line{\small\pfoot\ifnum\count0=\startpage
\copyright\ \gtp\hfill\else
\agt, Volume \thevolumenumber\ (\thevolumeyear)\hfill\fi}\vss}}
\else
\font=cmsl9 scaled 1050
\font\lnum=cmbx10 
\font=cmsl9 scaled 1050
\makeatletter
\def\@oddhead{{\small\lhead\ifnum\count0=\startpage ISSN 1472-2739 
(on-line) 1472-2747 (printed)\hfill {\lnum\number\count0}\else\ifodd\count0
\def\\{ }\ifx\theshorttitle\relax \thetitle \else\theshorttitle\fi\hfill
{\lnum\number\count0}\else\def\\{ and }{\lnum\number\count0}
\hfill\ifx\theshortauthors\relax 
\theauthors\else\theshortauthors\fi\fi\fi}}\def\@evenhead{\@oddhead}
\def\@oddfoot{\small\lfoot\ifnum\count0=\startpage\copyright\ \gtp\hfill\else
\agt, Volume \thevolumenumber\ (\thevolumeyear)\hfill\fi}
\def\@evenfoot{\@oddfoot}
\makeatother
\fi
\let\maketitlepage\makeagttitle
\let\makeshorttitle\maketitlepage
\let\maketitle\maketitlepage

\newwrite\gtoutfile
\long\gdef\makeheadfile{  
{\def\\{, }\def\s{ }
\immediate\openout\gtoutfile head.xxx
\immediate\write\gtoutfile{To: math@arxiv.org}
\immediate\write\gtoutfile{Subject: put OR rep NNNNN:ppppp}
\immediate\write\gtoutfile{--text follows this line--}
\immediate\write\gtoutfile{Proxy-for: \ifx\theasciiauthors\relax
\theauthors\else\theasciiauthors\fi\s<\ifx\theasciiemail\relax\theemail\else\theasciiemail\fi>}
\immediate\write\gtoutfile{\noexpand\\}
\immediate\write\gtoutfile{Authors: \ifx\theasciiauthors\relax
\theauthors\else\theasciiauthors\fi}
{\def\\{ }\immediate\write\gtoutfile{Title: \ifx\theasciititle\relax
\thetitle\else\theasciititle\fi}}
\immediate\write\gtoutfile{Subj-class: GT or SG, GR etc}
\immediate\write\gtoutfile{MSC-class: \theprimaryclass\ifx\thesecondaryclass\relax\else, \thesecondaryclass\fi}
\immediate\write\gtoutfile{Journal-ref: Algebr. Geom. Topol. \thevolumenumber\s
(\thevolumeyear) \startpage-\finishpage}
\immediate\write\gtoutfile{Comments: Published by Algebraic and
Geometric Topology at}
\immediate\write\gtoutfile{\s\s\s  http://www.maths.warwick.ac.uk/agt/AGTVol\thevolumenumber/agt-\thevolumenumber-\thepapernumber.abs.html}
\immediate\write\gtoutfile{\noexpand\\}
\immediate\write\gtoutfile{}
\ifx\theasciiabstract\relax
\immediate\write\gtoutfile{\theabstract}\else
\immediate\write\gtoutfile{\theasciiabstract}\fi
\immediate\write\gtoutfile{}
\immediate\write\gtoutfile{\noexpand\\}
\immediate\write\gtoutfile{}
\immediate\closeout\gtoutfile}}  

\def\maketitlepage{\makeagttitle\makeheadfile}
\let\makeshorttitle\maketitlepage
\let\maketitle\maketitlepage

\lognumber{46}
\volumenumber{2}
\volumeyear{2002}
\papernumber{46}
\pagenumbers{1179}{1195}
\received{21 June 2002}
\accepted{17 December 2002}
\published{27 December 2002}

\usepackage{amsmath,amssymb,epsf}

\let\nc\newcommand
\let\dmo\DeclareMathOperator
\let\nt\newtheorem

 \dmo{\mcg}{Mod}
 \dmo{\pmcg}{PMod} 
 \dmo{\homo}{H_1}
 \dmo{\aut}{Aut} 
 \dmo{\intnum}{i} 
 \dmo{\ain}{\hat{i}} 
 \dmo{\slz}{SL_2({\mathbb Z})} 
 \dmo{\relp}{gcd}

 \nc{\pic}[3]{\begin{figure}[ht!] \center{\leavevmode 
\epsfxsize #3 \epsfbox{#1.eps}} \caption{#2 \label{#1pic}} \end{figure}}

 \nc{\homs}{\ensuremath{\homo(S)}}
 \nc{\homsp}{\ensuremath{\homo(S')}} 
 \nc{\ra}{\rightarrow} 
 \nc{\bs}{\medskip} 
 
 \nc{\szfo}{\Sigma_{0,4}} 
 \nc{\bpf}{\begin{proof}} 
 \nc{\epf}{\end{proof}} 
 
 \nc{\ttt}[9]{\ensuremath{\left(\begin{array}{rrr} #1 & #2 & #3 \\ 
#4 & #5 & #6 \\ #7 & #8 & #9  \end{array}\right)}}

 \nc{\tbt}[4]{\ensuremath{\left(\begin{array}{rr} #1 & #2 \\
#3 & #4 \\ \end{array}\right)}}
 
 \nc{\tbo}[2]{\ensuremath{\left(\begin{array}{r} #1 \\ #2 \\ 
\end{array}\right)}}
 
 \nc{\thbo}[3]{\ensuremath{\left(\begin{array}{r} #1 \\ #2 \\
#3 \\ \end{array}\right)}}

 \nc{\tx}{\ttt{1}{-1}{-1}{0}{1}{0}{0}{0}{1}}
 \nc{\txi}{\ttt{1}{1}{1}{0}{1}{0}{0}{0}{1}}
 \nc{\ty}{\ttt{3}{-1}{-1}{2}{0}{-1}{2}{-1}{0}}
 \nc{\tyi}{\ttt{-1}{1}{1}{-2}{2}{1}{-2}{1}{2}}
 \nc{\typ}{\ttt{-1}{-1}{-1}{2}{2}{1}{2}{1}{2}}
 \nc{\typi}{\ttt{3}{1}{1}{-2}{0}{-1}{-2}{-1}{0}}
 \nc{\txty}{\ttt{-1}{0}{0}{2}{0}{-1}{2}{-1}{0}}
 \nc{\txtyp}{\ttt{-5}{-4}{-4}{2}{2}{1}{2}{1}{2}}
 \nc{\txtyi}{\ttt{3}{-2}{-2}{-2}{2}{1}{-2}{1}{2}}
 \nc{\txtypi}{\ttt{7}{2}{2}{-2}{0}{-1}{-2}{-1}{0}}

\nc{\subsec}{\subsection}

 \nt{lem}{{\bf Lemma}}{}
 \nt{thm}{{\bf Theorem}}{} 
  
 \nt{prop}{{\bf Proposition}}{} 
 \nt{form}{{\bf Formula}}{} 
 \nt{fact}{{\bf Fact}}{}

 \theoremstyle{definition}

 \nt{step}{Step \bs}

 \nc{\bstep}{\begin{step}}
 \nc{\estep}{\end{step}}

 \nc{\bl}{ \begin{list}{$\cdot$}{
 \setlength{\leftmargin}{.5in}
 \setlength{\rightmargin}{.5in}
 \setlength{\parsep}{0.5ex plus .2ex minus 0ex}
 \setlength{\itemsep}{0ex plus 0.2ex minus 0ex} 
 }
 }

 \nc{\el}{\end{list}}

\begin{document}

\title{A lantern lemma}

\authors{Dan Margalit}

\address{Department of Mathematics, 5734 S University Ave\\Chicago, 
IL 60637-1514, USA}                  

\email{margalit@math.uchicago.edu}                     

\url{http://www.math.uchicago.edu/margalit}                       

\begin{abstract}
We show that in the mapping class group of a surface any relation 
between Dehn twists of the form $T_x^jT_y^k=M$ ($M$ a multitwist) 
is the lantern relation, and any relation of the form 
$(T_xT_y)^k=M$ (where $T_x$ commutes with $M$) is the 2-chain 
relation. 
\end{abstract}

\asciiabstract{We show that in the mapping class group of a surface
any relation between Dehn twists of the form T_x^j T_y^k = M (M a
multitwist) is the lantern relation, and any relation of the form
(T_x T_y)^k = M (where T_x commutes with M) is the 2-chain relation.}




\maketitle 

\section{Introduction}

There is interesting interplay between the algebraic and
topological aspects of the mapping class group of a surface. One
instance is the algebraic characterization of certain topological 
relations between Dehn twists. For example, consider the 
following well-known relations (see, e.g. \cite[Chapter 4]{nvi}):

\bs

\begin{tabular}{l|l|l}
 & geometric relation & algebraic relation \\ \hline
reflexiveness & $a = b$ & $T_a = T_b$ \\
disjointness relation & $\intnum(a,b)=0$ & $T_aT_b=T_bT_a$ \\
braid relation & $\intnum(a,b)=1$ & $T_aT_bT_a=T_bT_aT_b$ \\
\end{tabular}

\bs

\noindent Here $a$ and $b$ are isotopy classes of simple closed
curves on a surface, $T_a$ and $T_b$ are the corresponding Dehn
twists, and $\intnum(a,b)$ is the geometric intersection number
between $a$ and $b$ (see below).

One can check directly that the given topological relations imply
the corresponding algebraic relations.  The algebraic relations 
{\em characterize} the topological relations in the sense that 
the algebraic relations imply the geometric ones. In other words, 
the algebraic relations {\em only} come from specific 
configurations of curves on the surface (see Section~\ref{basic}).

Ivanov-McCarthy give even more general statements \cite{im}: 

\bl

\item $T_a^j=T_b^k$ if and only if $a=b$ and $j=k$.

\item $T_a^jT_b^k=T_b^kT_a^j$ if and only if $\intnum(a,b)=0$.

\item $T_a^jT_b^kT_a^j=T_b^kT_a^jT_b^k$ if and only if $j=k=\pm 
1$ and $\intnum(a,b)=1$.

\el

McCarthy recently asked whether there was a similar
characterization of the {\em lantern relation}, which is a
relation between Dehn twists about curves which lie on a sphere
with four punctures:

\bs

\noindent {\bf Lantern relation}\qua $T_xT_yT_z =
T_{b_1}T_{b_2}T_{b_3}T_{b_4}$ where the curves are as in
Figure~\ref{lanternpic}.

\bs

Theorem~\ref{lantern} answers the question in the affirmative.

The lantern relation was discovered by Dehn \cite[Section 7g]{md},
and later by Johnson \cite[Section IV]{dlj}. Its significance
arises in part from the fact that it is one of very few relations
needed to give a finite presentation of the mapping class group
with the finite generating set of Humphries.

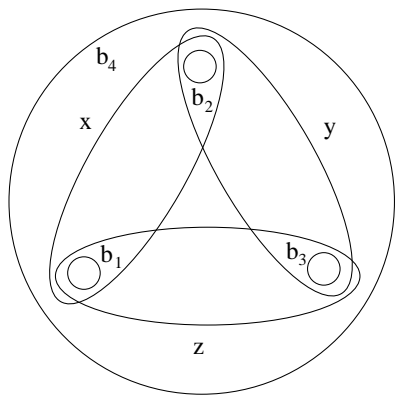
\pic{lantern}{The lantern configuration on a sphere with four
punctures}{2.05in}

For the statement of Theorem~\ref{lantern}, recall that the 
lantern relation can be written as $T_xT_y=M$, where $M$ is a {\em 
multitwist} (see Section~\ref{prelim}).

\begin{thm}[Lantern characterization]\label{lantern} Suppose
$T_x^jT_y^k=M$, where $M$ is a multitwist word and $j,k \in
{\mathbb Z}$, is a nontrivial relation between Dehn twists in
$\mcg(S)$.  Then the given relation (or its inverse) is the 
lantern relation; that is, $j=k=1$, a regular neighborhood of $x 
\cup y$ is a sphere with four boundary components, and 
$M=T_{b_1}T_{b_2}T_{b_3}T_{b_4}T_z^{-1}$, where the $b_i$ are the 
boundary components of the sphere, and $z$ is a curve on the 
sphere which has geometric intersection number 2 with both $x$ 
and $y$ (the sequence of curves $x$, $y$, $z$ should move 
clockwise around the punctured sphere as in 
Figure~\ref{lantern}). \end{thm}

\noindent In the theorem, the {\em inverse} of a relation 
$w_1=w_2$ between two words $w_1$ and $w_2$ is the equivalent 
relation $w_1^{-1}=w_2^{-1}$.

We also prove a similar theorem for the relation $(T_aT_b)^6=T_c$,
where $\intnum(a,b)=1$ and $c$ is the class of the boundary of a
regular neighborhood of $a \cup b$.

\begin{thm}[2-chain characterization]\label{2chain} Suppose
$M=(T_xT_y)^k$, where $M$ is a multitwist word and $k \in {\mathbb
Z}$, is a nontrivial relation between powers of Dehn twists in
$\mcg(S)$,  and $[M,T_x]=1$.  Then the given relation is the
2-chain relation---that is, $M=T_c^j$, where $c$ is the boundary 
of a neighborhood of $x \cup y$,  and $k=6j$\end{thm}

In Section~\ref{prelim} we prove the characterizations of the
disjointness relation and the braid relation, and introduce ideas
required for the proofs of our theorems.  Section~\ref{proof} is a
proof of Theorem~\ref{lantern} in the case $j=k=1$,
Section~\ref{genproof} generalizes to arbitrary $j$ and $k$,
Section~\ref{2cproof} contains the proof of Theorem~\ref{2chain},
Section~\ref{tech} contains supporting lemmas, and
Section~\ref{further} contains further questions related to this
work.

\paragraph{Acknowledgements} The author would like to thank John 
McCarthy for posing the problem, Benson Farb for relaying the 
problem and for discussing it in detail, and Feng Luo, Pallavi 
Dani, Angela Barnhill, and Kevin Wortman for helpful 
conversations and valuable input.  The author is also indebted to
the referee, Joan Birman, and Keiko Kawamuro for thoroughly 
reading a draft and making many suggestions and corrections.  
Hessam Hamidi-Tehrani, who has also been supportive, has obtained 
similar results \cite{hht}.

\section{Preliminaries}\label{prelim}

\subsec{Notation}

Let $S$ be an orientable surface.  We denote by $\mcg(S)$ the
mapping class group of $S$ (the group of orientation preserving
self-homeomorphisms of $S$, modulo isotopy).  When convenient, we
use the same notation for a curve on $S$, its isotopy class, and
its homology class.  Brackets around a curve will be used to
denote the homology class of the curve.

For two isotopy classes of simple closed curves $a$ and $b$ on
$S$, the {\em geometric intersection number} of $a$ and $b$,
denoted $\intnum(a,b)$, is the minimum number of intersection
points between representatives of the two classes.  By definition,
$\intnum(a,b)=\intnum(b,a)$.  The {\em algebraic intersection
number} of $a$ and $b$, denoted $\ain(a,b)$, is the sum of the
indices of the intersection points between any representatives of
$a$ and $b$, where an intersection point is of index $1$ when the
orientation of the intersection agrees with some given orientation
of the surface, and $-1$ otherwise.  Note that
$\ain(a,b)=-\ain(b,a)$.

If $a$ is an isotopy class of simple closed curves on $S$, we
denote by $T_a$ the mapping class of a Dehn twist about a
representative of $a$. As a matter of convention, Dehn twists will
be twists to the {\em left}.  Explicitly, if a neighborhood of a
representative of $a$ is an annulus $A$ parameterized (with
orientation) by $\{(r,\theta) \in {\mathbb R}^2:1 \leq r \leq 2, 0
\leq \theta \leq 2\pi\}$, then a representative of $T_a$ is the
diffeomorphism which is given by $(r,\theta) \mapsto (r, \theta +
(r-1)2\pi)$ on $A$ (in the given coordinates) and the identity
elsewhere.

A {\em multitwist} in $\mcg(S)$ is a product of Dehn twists
$\prod_{j=1}^n T_{a_j}^{e_j}$, where $\intnum(a_j,a_k)=0$ for any
$j$ and $k$ and $e_j \in {\mathbb Z}$.

The term {\em multitwist word} is used to describe a word in
$\mcg(S)$ consisting of Dehn twists about disjoint curves.  If
$M=\prod_{j=1}^n T_{a_j}^{e_j}$ is a multitwist word in $\mcg(S)$,
we can say that (for any $j$) the curve $a_j$ is {\em in M}.

\subsec{Formulas}

Ishida and Po\'enaru proved Formulas~\ref{relator} and~\ref{mt}, 
respectively, using elementary counting arguments \cite[Lemma 
2.1]{ai}\cite[Appendice Expos\'e 4]{flp}.  These inequalities are 
very useful in computations below.

\begin{form} \label{relator} Let $a$, $b$, and $c$ be any simple
closed curves on $S$, and let $n \in {\mathbb Z}$.  Then:
\[ |n|\intnum(a,b)\intnum(a,c)-\intnum(T_a^n(c),b) \leq
\intnum(b,c) \]
\end{form}

\begin{form} \label{mt} Let $M=\prod_{j=1}^n T_{a_j}^{e_j}$ be a
multitwist word with $e_j > 0$ for all $j$ (or $e_j < 0$ for all
$j$), and let $b$ and $c$ be arbitrary simple closed curves on
$S$. Then:
\[ |\intnum(M(c),b)-\Sigma_{j=1}^n e_j
\intnum(a_j,c)\intnum(a_j,b)| \leq \intnum(b,c) \]
\end{form}

As a special case of Formula~\ref{mt}, where $M=T_a^n$ and $b=c$,
we have:

\begin{form} \label{flp} Let $a$ and $b$ be any simple closed
curves on $S$.  Then:
\[ \intnum(T_a^{n}(b),b)=|n|\intnum(a,b)^2 \]
\end{form}

It will be essential in the proof of the Theorem~\ref{lantern} to
be able to compute the action of a product of Dehn twists on the
homology of a subsurface of $S$.  The well-known formula below is
the pertinent tool.

\begin{form} \label{hom} Let $a$ and $b$ be simple closed curves
on $S$, and $k$ an integer.  Then:
\[ [T_a^k(b)] = [b] + k\ain(b,a)[a] \]
where brackets denote equivalence classes in $\homs$. \end{form}

\subsec{Basic facts}

The following two facts are well-known and elementary
\cite[Corollary 4.1B, Lemma 4.1C]{nvi}.

\begin{fact} \label{unique} Let $a$ and $b$ be simple closed
curves on $S$.  If $T_a=T_b$, then $a$ is isotopic to
$b$.\end{fact}

\begin{fact}\label{conj} For $f \in \mcg(S)$ and $a$ any simple
closed curve on $S$, $fT_af^{-1}=T_{f(a)}$.\end{fact}

We now show that a Dehn twist about a given curve has a nontrivial
effect on every curve intersecting it.  This will be used, for
example, in the proof of Proposition~\ref{comm}.

\begin{fact}\label{effect}Let $a$ and $b$ be simple closed curves
on  on $S$.  If $\intnum(a,b) \neq 0$, then $T_a(b) \neq
b$.\end{fact}

\bpf  Using Formula~\ref{flp} we have
$\intnum(T_a(b),b)=\intnum(a,b)^2 \neq 0$.  On the other hand,
$\intnum(b,b)=0$.  Therefore $T_a(b) \neq b$.\epf

\subsec{Basic relation characterizations}\label{basic}

In the introduction, we stated characterizations of reflexiveness,
the disjointness relation, and the braid relation. The
characterization of reflexiveness is Fact~\ref{unique}. We present
the proofs of the latter two characterizations here for
completeness, and as a warmup for our main result.

\begin{prop}\label{comm}Let $a$ and $b$ be simple closed curves on
$S$.  If $T_aT_b=T_bT_a$ then $\intnum(a,b)=0$.\end{prop}

\bpf  Assume $T_aT_b=T_bT_a$.  Then, using Fact~\ref{conj}:
\begin{eqnarray*}
T_aT_b & = & T_bT_a \\
T_aT_bT_a^{-1} & = & T_b \\
T_{T_a(b)} & = & T_b 
\end{eqnarray*}
So $T_a(b)=b$ by Fact~\ref{unique}. By Fact~\ref{effect},
$\intnum(a,b)=0$. \epf

McCarthy proved the following characterization of the braid
relation in $\mcg(S)$ \cite[Lemma 4.3]{jdm}:

\begin{prop}\label{braid} Let $a$ and $b$ be non-isotopic simple
closed curves on $S$.  If $T_aT_bT_a=T_bT_aT_b$, then
$\intnum(a,b)=1.$\end{prop}

\bpf From the given algebraic relation and Fact~\ref{conj}, we
have:
\begin{eqnarray*}
T_aT_bT_a & = & T_bT_aT_b \\
(T_aT_b)T_a(T_aT_b)^{-1} & = & T_b \\
T_{T_aT_b(a)} & = & T_b \
\end{eqnarray*}
So $T_aT_b(a)=b$ by Fact~\ref{unique}.  Applying
Formula~\ref{flp}, we have:
\[ \intnum(a,b)^2 = \intnum(T_b(a),a) = \intnum(T_aT_b(a),a) =
\intnum(b,a) = \intnum(a,b) \]
Therefore $\intnum(a,b) \in \{0,1\}$.  If $\intnum(a,b)=0$, then
we have $T_a^2T_b=T_aT_b^2$, and hence $T_a=T_b$, i.e. $a$ is
isotopic to $b$, which contradicts the assumptions.\epf

Ivanov-McCarthy showed that actually the following more general
phenomena hold \cite[Theorems 3.14-3.15]{im}:

\begin{prop}\label{gref}Let $a$ and $b$ be simple closed curves on
$S$, and let $j$ and $k$ be nonzero integers. If $T_a^j=T_b^k$,
then $a$ is isotopic to $b$ and $j=k$.\end{prop}
\begin{prop}\label{gcomm}Let $a$ and $b$ be simple closed curves
on $S$, and let $j$ and $k$ be nonzero integers.  If
$T_a^jT_b^k=T_b^kT_a^j$ then $\intnum(a,b)=0$.\end{prop}
\begin{prop}\label{gbraid} Let $a$ and $b$ be non-isotopic simple
closed curves on $S$, and let $j$ and $k$ be nonzero integers.
If $T_a^jT_b^kT_a^j=T_b^kT_a^jT_b^k$, then
$\intnum(a,b)=1.$\end{prop}

\section{Proof of lantern characterization, $j=k=1$}\label{proof}

The idea is to build up, step by step, the lantern relation using
only the given algebraic information.  In particular, we show that
each of the following must be true for any algebraic relation
$T_xT_y=M$, where $M$ is a multitwist: $\intnum(x,y)>0$, $[M,T_x]
\neq 1$, there is a curve $z$ in the multitwist word $M$ with
$\intnum(x,z)>0$, $T_xT_y(z)=z$, $\intnum(x,z)=\intnum(y,z)$,
$\intnum(x,z)=\intnum(x,T_y(z))$, $\intnum(x,y)=2$, and
$\ain(x,y)=0$.  From this information, it will follow that the
given relation is the lantern relation.

\bstep\label{disjoint} $\intnum(x,y)>0$. \newline
If $\intnum(x,y)=0$, then $T_xT_y$ is a multitwist word, and so 
the multitwist word $M$ must also be $T_xT_y$ by 
Lemma~\ref{multitwists}, i.e. the equality between the words $M$ 
and $T_xT_y$ in $\mcg(S)$ is trivial.

\estep

\bstep\label{commute} $[M,T_x] \neq 1$. \newline
Assuming that $[M,T_x]=1$, we will arrive at a contradiction:
\[T_xT_yT_x^{-1}T_y^{-1}=MT_x^{-1}T_y^{-1}=T_x^{-1}MT_y^{-1} =
T_x^{-1}T_xT_yT_y^{-1}=1\]
So $T_xT_y=T_yT_x$, which implies $\intnum(x,y)=0$
(Proposition~\ref{gcomm}), contradicting Step~\ref{disjoint}.

\estep

\bstep\label{existz} There is a curve $z$ in the multitwist word 
$M$ with $\intnum(x,z)>0$. \newline
If $\intnum(x,z)=0$ for each curve $z$ in $M$, then $[M,T_x]=1$, 
which contradicts Step~\ref{commute}.  Therefore, there is a 
curve $z$ in $M$ which has nontrivial intersection with the curve 
$x$.

\estep

\bstep\label{fix} $T_xT_y(z)=z$. \newline
This is clear since $z$ is one of the curves in the multitwist 
word $M$: $T_xT_y(z)=M(z)=z$.

\estep 

\bstep\label{intnum} $\intnum(x,z)=\intnum(y,z)$. \newline 
Using Formula~\ref{flp}: $\intnum(T_y(z),z) = \intnum(y,z)^2$ and 
$\intnum(T_x^{-1}(z),z) = \intnum(x,z)^2$.  But since 
$T_x^{-1}(z) = T_y(z)$ (Step~\ref{fix}), all four expressions are 
the same, and so $ \intnum(y,z)^2 = \intnum(x,z)^2$.  Since 
geometric intersection number is a non-negative integer, we have 
$\intnum(x,z) = \intnum(y,z)$.

\estep 

\bstep\label{intnum2} $\intnum(y,z)=\intnum(x,T_y(z))$. \newline
Using $\intnum(T_y(z),z)=\intnum(y,z)^2$ (Formula~\ref{flp}) and 
$z=T_x T_y (z)$ (Step~\ref{fix}), we have: 
\[ \intnum(y,z)^2 = \intnum(T_y(z),z) = \intnum(z,T_y(z))= 
\intnum(T_x (T_y (z)),T_y(z)) = \intnum(x,T_y(z))^2 \]
So $\intnum(y,z) = \intnum(x,T_y(z))$.

\estep 

\bstep\label{intnum3} $\intnum(x,y)=2$. \newline 
Using Formula~\ref{relator}:
\[ \intnum(y,z)\intnum(y,x) - \intnum(x,T_y(z)) \leq \intnum(z,x) 
\]
But by Steps~\ref{intnum} and~\ref{intnum2}, 
$\intnum(y,z)=\intnum(x,z)=\intnum(x,T_y(z))$, so we can rewrite 
this as:
\[ \intnum(x,z)(\intnum(x,y)-2) \leq 0 \]
Since $\intnum(x,z) > 0$ (Step~\ref{existz}), this gives 
$\intnum(x,y) \in \{0,1,2\}$. The case $\intnum(x,y)=0$ is ruled 
out by Step~\ref{disjoint}.

We will now rule out $\intnum(x,y)=1$. In this case, a
neighborhood of $x \cup y$ on $S$ is a punctured torus $S'$.  We
will show that the induced action of $T_xT_y$ on $\homsp$, denoted
$(T_xT_y)_\star$, fails to fix any of the nontrivial elements of
$\homsp$; this contradicts the assumption that $T_xT_y$ is equal
to a multitwist in $\mcg(S)$ (Lemma~\ref{mixing}).

Using $\{[x],[y]\}$ as an ordered basis for $\homsp$, and
$\ain(x,y)=1$, Formula~\ref{hom} yields:
\[ (T_xT_y)_\star = (T_x)_\star(T_y)_\star = \tbt{1}{1}{0}{1} 
\tbt{1}{0}{-1}{1} = \tbt{0}{-1}{1}{1} \] 
This matrix does not have an eigenvalue of 1, so $(T_xT_y)_\star$ 
fixes no nontrivial element of $\homsp$.

Thus, $\intnum(x,y)=2$.

\estep 

\bstep\label{algint} $\ain(x,y)=0$. \newline
Since $\intnum(x,y)=2$, either $\ain(x,y)=0$ or $\ain(x,y)=\pm 
2$. We assume the latter and arrive at a contradiction.

Assuming $\ain(x,y)=\pm 2$ and $\intnum(x,y)=2$, a neighborhood of
$x \cup y$ (call it $S'$) is a genus one surface with two boundary
components (Figure~\ref{intxy2pic}).

As in Step~\ref{intnum3}, we will show that $(T_xT_y)_\star$ (the
induced action of $T_xT_y$ on $\homsp$) does not fix any
nontrivial, {\em nonperipheral} (see Lemma~\ref{mixing}) class in
$\homsp$. This will again contradict the assumption that $T_xT_y$
is equal to a multitwist in $\mcg(S)$ (Lemma~\ref{mixing}).

Let $x$, $v$, and $w$  be generators of $\homsp$ with
$\ain(x,v)=\ain(x,w)=1$, such that the two boundary components of
$S'$ are in the homology classes $v-w$ and $w-v$
(Figure~\ref{intxy2pic}).

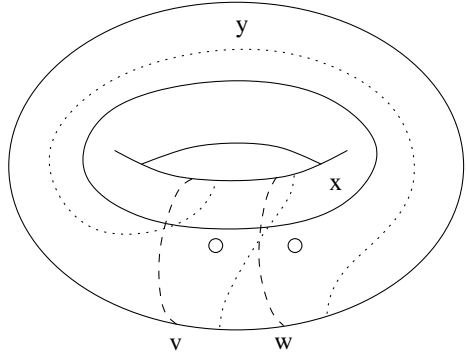
\pic{intxy2}{The picture for two simple closed curves with
algebraic intersection number~2}{2.4in}

Applying Formula~\ref{hom} and using $y=x+v+w$ (the case 
$\ain(x,y)=+2$), the action of $(T_xT_y)_\star = 
(T_x)_\star(T_y)_\star$ on $\homsp$ (with ordered basis 
$\{x,v,w\}$) is found to be:
\[ (T_xT_y)_\star = \tx \ty = \txty\]
In the case $\ain(x,y)=-2$, $y=x-v-w$ and the action is:
\[ (T_xT_y)_\star = \tx \typ = \txtyp\]
A basis for the fixed set of each of these linear operations is
$\{v-w\}$, which is the  homology class of a boundary component of
$S'$, i.e. the set of peripheral classes.

We have a contradiction, so $\ain(x,y)=0$.

\estep 

\bstep\label{endgame} The relation $T_xT_y=M$ is the lantern 
relation. \newline
Since $x$ and $y$ have geometric intersection number $2$ 
(Step~\ref{intnum3}) and algebraic intersection number $0$ 
(Step~\ref{algint}), a neighborhood of $x \cup y$ is a sphere 
with four boundary components $S'$. Let $M$ be the word 
\[ T_{b_1}T_{b_2}T_{b_3}T_{b_4}T_z^{-1}\]
where $b_i$ are the four boundary components of $S'$,
and $z$ is one of the two simple closed curves on $S'$ that hits 
each of $x$ and $y$ twice (the one pictured in 
Figure~\ref{lanternpic}), then it is well-known that $T_xT_y=M$ 
(To check this, draw any three arcs which cut $S'$ into a disk, 
and see that $T_xT_y$ and $M$ have the same effect on each of 
these arcs.  Then apply the Alexander lemma, which says that the 
mapping class group of a disk is trivial). By 
Lemma~\ref{multitwists}, $M$ is uniquely written as a product of 
twists about disjoint curves, and we are done.

\estep 

\section{Proof of general lantern 
characterization}\label{genproof}

\setcounter{step}{-1}

To show that any relation of the form $T_x^jT_y^k=M$ (where $M$ is
a multitwist word, $j,k \in {\mathbb Z}$) is the lantern relation,
we use the same program as in the proof for the case $j=k=1$ for
the first~\ref{gineq} steps.  Then, instead of homing in on
$\intnum(x,y)$, and $\ain(x,y)$, we show that $j=k=1$, which
leaves us in the case of Section~\ref{proof}.

\bstep\label{wlog} Assumptions on $j$ and $k$. \newline
We only consider ordered pairs of exponents $(j,k)$ in the set 
$\{(j,k):j>0,k>0\} \cup \{(j,k):j>0>k\}$ because $T_x^jT_y^k$ is 
equal to a multitwist word if and only if its inverse 
$T_y^{-k}T_x^{-j}$ is equal to a multitwist word.  Also, we can 
assume that both $j$ and $k$ are nonzero, because if at least one 
of them is zero, then $T_x^jT_y^k$ is a multitwist about one or 
no curves, and the relation $T_x^jT_y^k=M$ is trivial by 
Lemma~\ref{multitwists}.

\estep 

Steps~\ref{gdisjoint} through~\ref{gfix} are exactly the same as
for the case $j=k=1$, so we omit the proofs.

\bstep\label{gdisjoint} $\intnum(x,y)>0$.\estep 
\bstep\label{gcommute} $[M,T_x] \neq 1$. \estep 
\bstep\label{gexistz} There is a curve $z$ in the multitwist word 
$M$ with $\intnum(x,z)>0$. \estep
 \bstep\label{gfix} $T_x^jT_y^k(z)=z$. \estep

\bstep\label{gintnum} $|k|\intnum(y,z)^2=|j|\intnum(x,z)^2$. 
\newline
Using Formula~\ref{flp}: $\intnum(T_y^k(z),z) = 
|k|\intnum(y,z)^2$ and $\intnum(T_x^{-j}(z),z) = 
|j|\intnum(x,z)^2$.  Since $T_x^{-j}(z) = T_y^k(z)$ 
(Step~\ref{gfix}), all four expressions are equal, so we have
$|j|\intnum(x,z)^2 = |k|\intnum(y,z)^2$ and 
$\intnum(y,z)=\sqrt{|j/k|}\intnum(x,z)$.

\estep 

\bstep\label{gintnum2} $\intnum(x,T_y^k(z))=\intnum(x,z)$. 
\newline
Applying Step~\ref{gintnum}, Formula~\ref{flp}, 
Step~\ref{gfix}, and again Formula~\ref{flp}, we have: 
\[ |j|\intnum(x,z)^2 = |k|\intnum(y,z)^2 = \intnum(z,T_y^k(z))= 
\intnum(T_x^j(T_y^k(z)),T_y^k(z)) = |j|\intnum(x,T_y^k(z))^2 \]
So $|j|\intnum(x,z)^2=|j|\intnum(x,T_y^k(z))^2$, and further 
$\intnum(x,T_y^k(z))=\intnum(x,z)$.

\estep 

\bstep\label{gineq}$\intnum(x,y) \leq 2/\sqrt{|jk|}$. \newline
Using Formula~\ref{relator}: 
\[ |k|\intnum(y,z)\intnum(y,x) - 
\intnum(x,T_y^k(z)) \leq \intnum(z,x) \]
 But $\intnum(y,z)=\sqrt{|j/k|}\intnum(x,z)$ (Step~\ref{gintnum}) and 
$\intnum(x,T_y^k(z))=\intnum(x,z)$ (Step~\ref{gintnum2}), so we 
can rewrite this as: 
\[ \intnum(x,z)(\sqrt{|jk|}\intnum(x,y)-2) 
\leq 0 \]
Since $\intnum(x,z) > 0$ (Step~\ref{gexistz}), this gives 
$\intnum(x,y) \leq 2/\sqrt{|jk|}$.

\estep 

\bstep\label{gjk} $0 < |jk| \leq 4$. \newline
If $|jk|>4$, then the inequality of Step~\ref{gineq} says 
$\intnum(x,y) < 1$, which contradicts Step~\ref{gdisjoint}.  The 
inequality $|jk| > 0$ is part of Step~\ref{wlog}.

\estep 

\bstep\label{gjk1} $(j,k)=(1,\pm 1)$. \newline
If $(j,k) \neq (1,\pm 1)$, then $|jk|>1$ and Step~\ref{gineq} 
implies that $\intnum(x,y) < 2$.  This, coupled with 
$\intnum(x,y)>0$  (Step~\ref{gdisjoint}), gives $\intnum(x,y)=1$. 
In this case, a neighborhood of $x \cup y$ is a punctured torus 
$S'$, and $(T_x^jT_y^k)_\star$ acts on $\homsp$ (with basis 
elements represented by $x$ and $y$) via the matrix:
\[ (T_x^jT_y^k)_\star = \tbt{1}{1}{0}{1}^j\tbt{1}{0}{-1}{1}^k =
\tbt{1}{j}{0}{1} \tbt{1}{0}{-k}{1} = \tbt{1-jk}{j}{-k}{1}\]
which has eigenvalues:
\[ e(j,k) = \frac{(2-jk) \pm \sqrt{(jk)^2-4jk}}{2} \] 
By Lemma~\ref{mixing}, since $T_x^jT_y^k$ is equal to a 
multitwist supported on $S'$, $(T_x^jT_y^k)_\star$ must have a 
fixed point on $\homsp$, so it must have an eigenvalue of $1$.  
We will show, however, that $e(j,k)$ does not equal 1 for $1 < 
|jk| \leq 4$, and so $|jk|=1$.

By the standing assumption that either $j$ and $k$ are positive 
or $j>0>k$ (Step~\ref{wlog}), and the fact that $e(j,k) = e(k,j) 
= e(1,jk)$, it suffices to check $e(1,jk)$ for $2 \leq |jk| \leq 
4$.  Using the formula, we have $e(1,4)=-1$, $e(1,3)=(-1 \pm 
\sqrt{3}i)/2$, $e(1,2)=\pm i$, $e(1,-2)= 2 \pm \sqrt{3}$, 
$e(1,-3)=(5 \pm \sqrt{21})/2$, and $e(1,-4)=3 \pm 2 \sqrt{2}$.

\estep 

\bstep\label{gjk12} $(j,k)=(1,1)$. \newline
By Step~\ref{gjk1}, the only possibilities left for $(j,k)$ are
$(1,1)$ and $(1,-1)$.   Our goal now is to show that $(j,k)=
(1,-1)$ leads to a contradiction.  Step~\ref{gineq} implies that
$\intnum(x,y) \leq 2$ in this case.

As in Step~\ref{gjk1}, we know $\intnum(x,y) \neq 1$ because
$e(1,-1)=(3 \pm \sqrt{5})/2$. In particular 
$(T_x^1T_y^{-1})_\star$ does not have an eigenvalue of $1$, 
contradicting Lemma~\ref{mixing}.

We can also check that $\intnum(x,y) \neq 2$.  There are three
subcases: $\ain(x,y)=2$, $\ain(x,y)=-2$, and $\ain(x,y)=0$.

For $\ain(x,y)=2$, we can compute $(T_xT_y^{-1})_\star = 
(T_x)_\star(T_y)_\star^{-1}$ as follows (using the ordered basis 
$\{x,v,w\}$ as in Section~\ref{proof}, Step~\ref{algint}):
\[  (T_xT_y^{-1})_\star = \tx \tyi = \txtyi   \]
And for $\ain(x,y)=-2$, we have:
\[  (T_xT_y^{-1})_\star = \tx \typi = \txtypi   \]
The only fixed points of the above two matrices are peripheral
classes (multiples of $v-w$).  By Lemma~\ref{mixing}, both of
these cases are impossibilities.

The final subcase for $(j,k)=(1,-1)$ and $\intnum(x,y)=2$ is
$\ain(x,y)=0$. In this situation, a regular neighborhood of $x
\cup y$ is a sphere with four punctures $S'$. Since $\homsp$ 
contains only peripheral elements, Lemma~\ref{mixing} does not 
apply.  We employ a similar idea, with curve classes playing the 
role of homology classes. In particular, we will show that 
$T_xT_y^{-1}$ is {\em irreducible} on $S'$ (i.e. it does not fix 
any nontrivial isotopy class of simple closed curves on $S'$). By 
Lemma~\ref{irred}, this is a contradiction.

It is well known that the isotopy classes of simple closed curves 
on $S'$ are in one-to-one correspondence with the set $\{(p,q): 
\relp(p,q)=1\}/\sim$, where $(p,q)\sim(-p,-q)$, that $\pmcg(S')$ 
(the subgroup which preserves punctures) is isomorphic to a 
finite-index subgroup of $\slz$ with a matrix $A$ acting on a 
$(p,q)$ curve by matrix multiplication, and that a Dehn twist 
about the $(1,0)$ curve is given by the matrix $((1,2),(0,1))$ 
\cite[Section 3]{ynm}. Therefore:
\[ T_xT_y^{-1} = \tbt{1}{2}{0}{1} \tbt{1}{0}{-2}{1}^{-1}=
\tbt{1}{2}{0}{1} \tbt{1}{0}{2}{1} = \tbt{5}{2}{2}{1} \]
This matrix does not fix any $(p,q)$ (since it does not have an
eigenvalue of $\pm 1$). In other words, the mapping class is
irreducible, and by Lemma~\ref{irred} this contradicts the
assumption that $T_xT_y^{-1}$ is equal to a multitwist.

\estep 

\bstep\label{gendgame} The relation $T_x^jT_y^k=M$ is the lantern 
relation. \newline
We have eliminated all possibilities for the exponents except
$j=k=1$.  By Section~\ref{proof}, the given relation (or its 
inverse) is the lantern relation.

\estep 

\section{Proof of 2-chain characterization}\label{2cproof}

Theorem~\ref{2chain} follows from a result proven by Ishida and
Hamidi-Tehrani \cite[Theorem 1.2]{ai} \cite{hht}:

\bs \noindent {\bf Theorem}\qua {\sl If $\intnum(x,y) \geq 2$, then
there are no relations between $T_x$ and $T_y$.} \bs

If we have the conditions of the theorem: $(T_xT_y)^k=M$, where
$M$ is a multitwist word, and $[T_x,M]=1$.  Then:
\[ (T_xT_y)^kT_x = MT_x=T_xM=T_x(T_xT_y)^k \]
which is a relation between $T_x$ and $T_y$,  assuming $|k|>1$ (by
Section~\ref{proof}, Step~\ref{commute} there are no relations
with $|k|=1$ and  $[T_x,M]=1$) . Therefore, $\intnum(x,y) \in
\{0,1\}$. We can rule out $\intnum(x,y)=0$, because then the
relation $T_xT_y=M$ is trivial by Lemma~\ref{multitwists}.  Thus,
$\intnum(x,y)=1$, and a neighborhood of $x \cup y$ is a punctured
torus $S'$.

As in Section~\ref{proof}, Step~\ref{intnum3}, we consider the
action of $(T_xT_y)^k$ on $\homsp$ with generators represented by
$x$ and $y$.  The first $6$ powers of $(T_xT_y)_\star$ are
\[ \tbt{0}{1}{-1}{1}, \tbt{-1}{1}{-1}{0}, \tbt{-1}{0}{0}{-1}, 
\tbt{0}{-1}{1}{-1},\tbt{1}{-1}{1}{0}, \tbt{1}{0}{0}{1} \]
The first five of these matrices fix no nontrivial vector.  
Hence, by Lemma~\ref{mixing}, $(T_xT_y)^k$ cannot equal a 
multitwist in $\mcg(S)$ for $k$ not a multiple of $6$.  When 
$k=6j$ for some integer $j$, then it is well-known that 
$(T_xT_y)^k=T_c^j$ where $c$ is the boundary component of $S'$ 
\cite[Lemma 4.1G]{nvi}.  One can check this relation by using the 
Alexander lemma, as in Section~\ref{proof}, Step~\ref{endgame}. 
By Lemma~\ref{multitwists}, the multitwist word $M$ is unique, 
and we are done.

\section{Technical lemmas} \label{tech}

Lemma~\ref{multitwists} uses some new terminology: An {\em 
essential reduction class} of $f \in \mcg(S)$ is a class of 
simple closed curves $\alpha$ such that $f(\alpha)=\alpha$, and 
if $\intnum(\alpha,\gamma) \neq 0$ then $f^n(\gamma)\neq\gamma$ 
for any $n \in {\mathbb N}$.  The {\em canonical reduction 
system} for $f \in \mcg(S)$ is the set of essential reduction 
classes of $f$. Lemma~\ref{multitwists} is really a special case 
of the theorem of Birman-Lubotzky-McCarthy which states that 
canonical reduction systems are unique.

\begin{lem}\label{multitwists}
Suppose $M=\prod_{j=1}^m T_{x_j}^{e_j}$ and $N=\prod_{j=1}^n
T_{y_j}^{f_j}$ are multitwist words in $\mcg(S)$.  If $M=N$ in
$\mcg(S)$, then $m=n$ and $\{(x_j,e_j)\}= \{(y_j,f_j)\}$.\end{lem}

\bpf Since $M$ and $N$ are multitwist words, 
$\intnum(x_j,x_k)=\intnum(y_j,y_k)=0$, and so $M(x_i)=x_i$ and 
$N(y_j)=y_j$ for all $i$ and $j$.  It then follows from the work 
of Birman-Lubotzky-McCarthy that both $\{x_j\}$ and $\{y_j\}$ are 
canonical reduction systems for $M=N$ \cite[Lemma 2.5]{blm}, and 
hence the sets are the same by uniqueness of such systems 
\cite[Theorem C]{blm}.  It then follows that the exponents are 
the same: consider the surface obtained by cutting $S$ along 
$\{x_k\}_{k \neq j}$;  the mapping class induced by $M$ on this 
surface is $T_{x_j}^{e_j}=T_{x_j}^{f_j}$ (assuming $x_j=y_j$), 
and no two different powers of a Dehn twist are the same element 
(Proposition~\ref{gcomm}), so $e_j=f_j$. \epf

For Lemma~\ref{mixing}, a {\em peripheral homology class} $\alpha
\in \homsp$ on a subsurface $S' \subset S$ is one which is 
contained in the subgroup of $\homsp$ generated by the classes of
components of $\partial S'$.  For $f \in \mcg(S)$, $f_\star$ 
denotes the induced action of $f$ on homology.

\begin{lem}\label{mixing} Suppose $M \in \mcg(S)$
is a multitwist with support on a subsurface $S'$, and that there
is a nontrivial and nonperipheral element of $\homsp$.  Then there
is a nontrivial and nonperipheral $\alpha \in \homsp$ with
$M_\star(\alpha)=\alpha$.\end{lem}

\bpf Since $M$ has its support on $S'$, it must be of the form:
\[ M = \prod_{i=1}^m T_{a_i}^{e_i} \prod_{j=1}^n T_{b_j}^{f_j}
\prod_{k=1}^p T_{c_k}^{g_k} \]
where the $a_i$ represent the trivial class in $\homsp$, the 
$b_j$ represent peripheral homology classes in $\homsp$, and the 
$c_k$ represent nontrivial, nonperipheral classes in $\homsp$.  
If $p$ is nonzero, i.e. $M$ consists of at least one twist about 
a representative of a nontrivial, nonperipheral homology class, 
then $M_\star([c_k])=[c_k]$ for any $k$ since $M(c_k)=c_k$, and 
we are done. Otherwise, if $p=0$, let $s$ be a simple closed 
curve on $S'$ representing any nontrivial, nonperipheral class in 
$\homsp$.  Then $\ain([s],[a_i])=0$ (the $[a_i]$ can be 
represented by the trivial curve class) and $\ain([s],[b_i])=0$ 
(the $[b_i]$ can be represented by boundary curves), so 
$M_\star([s])=[s]$ by Formula~\ref{hom}.
\epf

Recall that an irreducible mapping class is one which fixes no
isotopy class of curves.  Lemma~\ref{irred} states that a
multitwist in $\mcg(S)$ cannot restrict to an irreducible mapping
class on a subsurface of $S$.

\begin{lem}\label{irred} Suppose $M \in \mcg(S)$ is a multitwist
with support on a subsurface $S'$, and that there is a nontrivial
(not homotopic to a point or a boundary component) isotopy class
of curves on $S'$.  Then there is a nontrivial isotopy class of
curves on $S'$ which is fixed by $M$.\end{lem}

\bpf

Since $M$ has support on $S'$, it is of the form:
\[ M = \prod_{i=1}^m T_{a_i}^{e_i} \prod_{j=1}^n T_{b_j}^{f_j} \]
where the $a_i$ are boundary components of $S'$ and the $b_j$ are
nontrivial curve classes on $S'$.  If $n \neq 0$, then
$M(b_j)=b_j$ for any $1 \leq j \leq n$. If $n = 0$, then
$M(\alpha)=\alpha$ for any nontrivial curve class on $S'$.
\epf

\section{Questions}\label{further}

\paragraph{Powers of $T_xT_y$} This paper gives a partial
classification of relations of the form $(T_xT_y)^k=M$, where $M$ 
is a multitwist word. If $k=1$, then it is the lantern relation. 
If $k \neq 1$ and $[T_x,M]=1$, then it is the 2-chain relation. 
The author is unaware of relations where $k \neq 1$ and $[T_x,M] 
\neq 1$. One way to generalize this is to consider relations of 
the form $W(T_x,T_y)=M$, where $W(T_x,T_y)$ is any word in $T_x$ 
and $T_y$.\footnote{Hamidi-Tehrani has successfully addressed 
this question \cite{hht}.}

\paragraph{Noncommutativity} The results of this paper rely
heavily on the assumption that certain mapping classes are
multitwists.  This is a strong assumption, as multitwists a priori
consist of {\em disjoint} curves.  A more general problem is to
classify all relations of the form $T_xT_yT_z=T_aT_bT_cT_d$, with
no hypotheses of commutativity or disjointness.

\paragraph{Multiple lanterns} A natural question to ask is under
what assumptions is $XY=M$ ($M$ a multitwist word) the lantern
relation for arbitrary mapping classes. This is certainly not true
for any $X$ and $Y$. For example, there are {\em multiple
lanterns}: Let $X=T_{x_1}T_{x_2}$ and $Y=T_{y_1}T_{y_2}$, where
$T_{x_1}T_{y_1}=M_1$ and $T_{x_2}T_{y_2}=M_2$ are lantern
relations.  Then $XY$ is a multitwist if $[M_1,M_2]=1$. If the two
lanterns have the same boundary components, then $M=M_1^2=M_2^2$.

\paragraph{Chain relations} There is a canonical relation for any
$n$-chain of curves on $S$ (a sequence of curves $\{a_1, \dots,
a_n\}$ with $\intnum(a_j,a_k)=1$ for $k=j \pm 1$ and
$\intnum(a_j,a_k)=0$ otherwise).  When $n$ is odd, the boundary of
a neighborhood of the $n$-chain consists of two curves $d_1$ and
$d_2$, and we have the relation: 
\[ (T_{a_1} \dots T_{a_n})^{n+1} = T_{d_1}T_{d_2} \]
and when $n$ is even, a neighborhood of the $n$-chain consists of 
one curve $d_1$, and we have:
\[ (T_{a_1} \dots T_{a_n})^{2n+2} = 
T_{d_1} \]
One can ask how well these relations can be characterized. Note 
that Theorem~\ref{2chain} is the special case $n=2$, and that the 
case of $n=1$ is Fact~\ref{unique}.

\Addresses\recd


\begin{thebibliography}

\bibitem{flp}
\emph{Travaux de {T}hurston sur les surfaces}, Soci\'et\'e Math\'ematique de
  France, Paris (1991), s\'eminaire Orsay, Reprint of {\it Travaux de Thurston
  sur les surfaces}, Soc.\ Math.\ France, Paris, 1979, Ast\'erisque No. 66-67
  (1991)

\bibitem{blm}
\textbf{Joan~S Birman}, \textbf{Alex Lubotzky}, \textbf{John McCarthy},
  \emph{Abelian and solvable subgroups of the mapping class groups}, Duke Math.
  J. 50 (1983) 1107--1120

\bibitem{md}
\textbf{Max Dehn}, \emph{Papers on group theory and topology}, Springer-Verlag,
  New York (1987)

\bibitem{hht}
\textbf{Hessam Hamidi-Tehrani}, \emph{{Groups generated by postitive
  multi-twists and the fake lantern problem}}, Algebr. Geom. Topol. 2
(2002) 1155--1178, {\sl immediately preceding this paper}

\bibitem{ai}
\textbf{Atsushi Ishida}, \emph{The structure of subgroup of mapping class
  groups generated by two {D}ehn twists}, Proc. Japan Acad. Ser. A Math. Sci.
  72 (1996) 240--241

\bibitem{nvi}
\textbf{Nikolai~V Ivanov}, \emph{Mapping class groups}, from: ``Handbook of
  geometric topology'', North-Holland, Amsterdam (2002)  523--633

\bibitem{im}
\textbf{Nikolai~V Ivanov}, \textbf{John~D McCarthy}, \emph{On injective
  homomorphisms between {T}eichm\"uller modular groups. {I}}, Invent. Math. 135
  (1999) 425--486

\bibitem{dlj}
\textbf{Dennis~L Johnson}, \emph{Homeomorphisms of a surface which act
  trivially on homology}, Proc. Amer. Math. Soc. 75 (1979) 119--125

\bibitem{jdm}
\textbf{John~D McCarthy}, \emph{Automorphisms of surface mapping class groups.
  {A} recent theorem of {N}. {I}vanov}, Invent. Math. 84 (1986) 49--71

\bibitem{ynm}
\textbf{Yair~N Minsky}, \emph{A geometric approach to the complex of curves on
  a surface}, from: ``Topology and Teichm\"uller spaces (Katinkulta, 1995)'',
  World Sci. Publishing, River Edge, NJ (1996)  149--158

\end{thebibliography}
\end{document}